\documentclass[12pt]{article}
\usepackage{amsmath}
\usepackage{amsfonts}
\usepackage{a4wide}

\usepackage{amssymb}
\usepackage{latexsym}
\usepackage{euscript}

\def\al{w^{*}}
\def\<{\langle}
\def\>{\rangle}
\def\ip{\<\ , \ \>}

\def\E{\mathcal E}

\def\qed{\hfill \square}

\def\Mo{\M_0}
\def\So{\mathbf{S}_0}
\def\Sk{\mathbf{S}}
\newcommand{\Z}{\mathbf{Z}}
\newcommand{\A}{\mathbf A}
\newcommand{\Af}{\A_f}
\newcommand{\C}{\mathbf C}
\newcommand{\F}{\mathbf F}
\newcommand{\Fbar}{\overline{\F}}
\newcommand{\Fp}{\F_p}

\newcommand{\Kbar}{\overline{K}}
\newcommand{\PP}{\mathbf P}
\newcommand{\Q}{\mathbf{Q}}

\newcommand{\Qp}{\Q_p}
\newcommand{\M}{\mathbf{M}}
\newcommand{\R}{\mathbf{R}}
\newcommand{\Zhat}{\widehat{\Z}}
\newcommand{\Zp}{\Z_p}
\DeclareMathOperator{\an}{an}
\DeclareMathOperator{\GL}{GL}
\DeclareMathOperator{\SL}{SL}
\DeclareMathOperator{\Spec}{Spec}

\title{Slopes of overconvergent 2-adic modular forms.}
\author{Kevin Buzzard \and Frank Calegari\footnote{Supported in part by the American Institute of Mathematics.}}

\begin{document}
\maketitle

\newtheorem{theorem}{Theorem}
\newtheorem{df}{Definition}
\newtheorem{lemma}{Lemma}
\newtheorem{cor}{Corollary}
\newtheorem{conj}{Conjecture}
\newtheorem{proposition}{Proposition}

\section{Introduction.}

Let~$p$ be a prime,
and let~$N$ be a positive integer coprime to~$p$.
Let $M_k(\Gamma_1(N);\Qp)$ denote the weight~$k$ modular forms of
level $\Gamma_1(N)$ defined over $\Qp$. 
In recent years, work of Coleman and others
(for example 
\cite{coleman:families},\cite{coleman2},\cite{coleman3},\cite{coleman4},\cite{eigencurve})
has shown that a very
profitable way of studying this finite-dimensional $\Qp$-vector space
is to choose a small positive rational number~$r$ and then to embed
$M_k(\Gamma_1(N);\Qp)$ into a (typically infinite-dimensional) $p$-adic Banach
space $\M_k(\Gamma_1(N);\Qp;p^{-r})$ of $p^{-r}$-overconvergent $p$-adic
modular forms,
that is, sections of $\omega^{\otimes k}$ on the affinoid subdomain
of $X_1(N)$ obtained by removing certain open discs of radius $p^{-r}$
above each supersingular point in characteristic~$p$ (at least if $N\geq5$;
see the appendix for how to deal with the cases $N\leq 4$).
The space $\M_k(\Gamma_1(N);\Q_p;p^{-r})$, for $0<r<p/(p+1)$,
comes equipped with canonical continuous Hecke operators,
and one of them, namely the operator $U:=U_p$,
has the property of being compact; in particular $U$ has a spectral
theory. Coleman exploited this theory in~\cite{coleman:families}
to prove weak versions
of conjectures of Gouv\^ea and Mazur on families of modular forms.

One of us (K.B.) has made, in many cases, considerably more precise
conjectures (\cite{buzzard:slopeconjectures}) than those of Gouv\^ea
and Mazur, predicting the slopes of $U$,
that is, the valuations of all the non-zero eigenvalues of $U$.
These conjectures are very explicit, and display a hitherto
unexpected regularity. However, they have the disadvantage of being rather
inelegant. See also the forthcoming PhD thesis~\cite{Herrick} of
Graham Herrick, who has, perhaps, more conceptual conjectures
about these slopes.

We present here a very concrete conjecture in the case $N=1$ and $p=2$,
which presumably agrees with the conjectures in
\cite{buzzard:slopeconjectures} but which has the advantage of
being much easier to understand and compute. Let $S_k:=S_k(\Gamma_0(1),\Q)$
 denote
the level~1 cusp forms of weight~$k$. If~$F(X)$ is a polynomial
with rational coefficients then by its 2-adic Newton polygon we mean
its Newton polygon when considered as a polynomial with 2-adic
coefficients.

\begin{conj} \label{conj:classical} Let $k \ge 12$ be even, and
let $m = \dim S_k$.
Then the 2-adic Newton Polygon of \mbox{$\mathrm{det}(1-X T_2)$} 
on $S_k$ equals the 2-adic Newton Polygon of  
$$1 + \sum_{n=1}^{m}
X^n \prod_{j=1}^n \frac{2^{2j} 
(k-8j)! (k-8j-3)! (k-12j-2)}{(k-12j)! (k-6j-1)!}.$$
\end{conj}
This conjecture can be verified numerically, and we have verified
it for all $k\leq2048$. Control theorems of Coleman (\cite{coleman2}) imply
that complete knowledge of all slopes of classical cusp forms at level~2
for all weights is equivalent
to complete knowledge of all slopes of finite overconvergent
tame level~1 cusp forms for all integer weights. In fact, it is a little
tedious but completely elementary to show that we may
reformulate
Conjecture~\ref{conj:classical} as follows. Let
$\Sk_k:=\Sk_k(\Gamma_0(1);\Q_2;2^{-1/2})$
denote the $2^{-1/2}$-overconvergent forms of weight~$k$ and
tame level~1.

\begin{conj} \label{conj:nonclassical} Let $k\leq0$ be an integer.
Then the Newton Polygon of \mbox{$\mathrm{det}(1-X U)$}
on $\Sk_k$ is
the Newton Polygon of
$$1 + \sum_{n=1}^{\infty}
X^n \prod_{j=1}^n \frac{2^{2j} (-k+2+12j)! (-k+6j)!}{(-k+2+8j)!(-k-2+8j)!(-k-12j)}.$$
\end{conj}

One can also find a form of this conjecture that makes sense if $k>0$,
for example if one reformulates the factorials as Gamma-functions and then
is careful to make precise what is happening at poles. We leave this
reformulation to the reader.

As evidence for this conjecture, we have the following:

\begin{theorem}\label{temp} Conjecture~\ref{conj:nonclassical} is
true when $k = 0$.
\end{theorem}

If $x$
is a non-zero rational number, then by its \emph{slope} we mean
its 2-adic valuation $v_2(x)$. 
It is easy to check that for $n\geq0$ an integer
we have $v_2\left((2n)!\right)=n+v_2(n!)$, and it
follows from this that
$$v_2\left( \frac{2^{2j} (12j + 2)! (6j)!}{(8j+2)!(8j-2)!(12j)}\right)
=  1 + 2 v_2 \left( \frac{(3j)!}{j!} \right).$$
Similarly one can check that if $D$ denotes the infinite diagonal
matrix whose $(j,j)$th entry, $j\geq1$, is given by
$$d_{j,j} = \frac{2^{4j + 1} (3j)!^2 j!^2}{3 \cdot (2j)!^4},$$
then $v_2(d_{j,j})=1 + 2 v_2 \left( \frac{(3j)!}{j!} \right).$
Hence Theorem~\ref{temp} above is equivalent to

\begin{theorem} \label{theorem:weight0} The Newton polygons
of $\det(1-XU)$ on $\So$ and $\det(1-XD)$ coincide.
\end{theorem}

This is the form of the theorem that we shall actually prove.

Note that the sequence $v_2((3j)!/j!)$ is strictly increasing;
this follows from the fact 
$$\frac{(3j+3)!/(j+1)!}{(3j)!/j!} = 3(3j+2)(3j+1)$$
is even for all $j$. We deduce

\begin{cor}\label{distinctslopes} 
 Let $|\lambda_1| \ge |\lambda_2|
 \ge |\lambda_3| \ge \ldots$
be the non-zero eigenvalues \emph{(}with multiplicities\emph{)}
 of $U$ on $\So$.  Then the slope of $\lambda_n$ is given by the following
formula:
$$v_2(\lambda_n) = 1 + 2 \cdot
v_2 \left( \frac{(3n)!}{(n)!}\right).$$
\end{cor}

In particular, the slopes are all distinct, and are all positive
odd integers. We also have

\begin{cor}
Let
$f=q+\ldots \in \So\widehat{\otimes}\C_2$
be a normalised finite slope overconvergent eigenform.
Then the coefficients of $f$ are all in~$\Q_2$. 
\end{cor}

\begin{Proof} We use only that the slopes of the non-zero eigenvalues
of~$U$ are distinct. If $\lambda$ denotes
the eigenvalue
of~$U$ on~$f$ then $1-\lambda^{-1}U$ is not invertible
on $\So\widehat{\otimes}\C_2$ and hence
by Proposition~11 of~\cite{serre} we see that~$\lambda^{-1}$ is a zero
of the characteristic power series~$P(T)$ of~$U$ acting on~$\So$.
Note that $P(T)\in\Q_2[[T]]$. Choosing some big affinoid disc
containing~$\lambda^{-1}$ and applying the Weierstrass preparation theorem
to~$P(T)$ shows that $\lambda^{-1}$ is a root of a polynomial with
coefficients in~$\Q_2$.
Hence $\lambda\in\overline{\Q}_2$. Now all the Galois conjugates of~$\lambda$
have the same valuation and are also roots of the characteristic power
series of~$U$; hence, by Corollary~\ref{distinctslopes}, $\lambda\in\Q_2$.
Finally by Proposition~12
of~\cite{serre} the subspace of~$\So$ where~$U$ acts as
multiplication by~$\lambda$ is one-dimensional over $\Q_2$ and so
all the eigenvalues of all the other Hecke operators are also in $\Q_2$.
$\qed$
\end{Proof}

\medskip

{\bf Remark}. This Corollary provides some evidence
towards Question $4.3$ of \cite{buzzard:slopeconjectures}. See
also Corollary~1.2 of~\cite{kilford:thesis}. 

\

Note that for $p=2$ and $N=1$, the map $\theta$ defines an isomorphism
$\theta: \So \rightarrow \Sk_2$ such that $U \theta = 2 \theta U$.
Thus, the slopes in weight two are precisely each of the slopes
in weight zero, plus one.
We have also proved Conjecture~\ref{conj:nonclassical}
for $k=-12$ using similar methods, although
the combinatorics are too painful to write here, and the arguments
do not seem to generalise to all $k$.

Lawren Smithline was perhaps
the first person to observe
that there was some structure in the slopes of overconvergent modular
forms of small level; his results~(\cite{smithline:thesis})
were primarily for the prime~$p=3$
but some of the techniques used in this paper for studying the explicit
matrix representing~$U$ were inspired by his ideas. As far as we know,
the first people to get explicit results pinning down all overconvergent
slopes at a given weight were L.~Kilford~(\cite{kilford:thesis})
and D.~Jacobs~(\cite{jacobs:thesis}) but their results differ in two
respects from ours: firstly, they consider points nearer the
boundary of weight space, and secondly the slopes at the weights they
consider have a much simpler pattern---they form an arithmetic
progression. 

We would like to thank Robert Coleman and Matthew
Emerton for several helpful
discussions. This genesis of this paper was a project
at the
\oldstylenums{2001} 
Arizona Winter School.

\section{Weight Zero}

The curve $X_0(2)$ has genus $0$. A natural choice of uniformiser
is given by the following function (Hauptmodul):
$$f(\tau) = \Delta(2 \tau)/\Delta(\tau)$$
In the sequel, we shall write this function simply
as $f$. There is
a product formula for $f$:
$$f = q \prod_{n=1}^{\infty} (1+q^n)^{24} = 
q \prod_{n=1}^{\infty} \frac{1}{(1-q^{2n-1})^{24}} =  q + 24 q^2 + 300 q^3 + 
2624 q^4 + 18126 q^5 + 105504 q^6 + \ldots$$
which follows immediately
 from the usual product formula for $\Delta$.
For $p = 2$, 
$f$ is overconvergent of weight $0$ and level $1$. Furthermore, if $g=2^6f$
then the set
$\{1,g,g^2,g^3,\ldots\}$ is a Banach basis for the space
$\Mo:=\M_0(\SL_2(\Z);\Q_2;2^{-1/2})$. It seems a little difficult
to extract these concrete statements from the literature and so we sketch
a proof of this in the appendix, and note that these ideas should
easily be adaptable to cover other cases where one might want to do
explicit computations. The reader who is happy to accept that
the formal Banach space with basis $\{1,g,g^2,g^3,\ldots\}$ is
some kind of $p$-adic space of modular forms might well want
to avoid these technical details.

To determine the spectral theory of $U$ on $\Mo$, we shall
explicitly compute a matrix for $U$ acting on $\Mo$. These
calculations are much in the  spirit  of classical
congruences for coefficients of modular
functions such as $j$, see for example
Watson \cite{Watson}, or Atkin and O'Brien \cite{At}. In this optic, the
\emph{a fortiori} presence of a spectral theory greatly
simplifies matters. Concretely then, our task is to compute $U(f^k)$
as a power series in $f$, for all $k\geq0$.

\begin{lemma} \label{lemma:UF} The following identities
are satisfied:
$$U(f) = 24 f + 2^{11} f^2.$$
$$f\bigg(\frac{\tau}{2}\bigg) f\bigg(\frac{\tau+1}{2}\bigg) = - f(\tau).$$
\end{lemma}

\begin{Proof} The operator $U$ preserves the space
of functions on $X_0(2)$. Furthermore $U(f)$, considered
as a map $X_0(2)\to\PP^1$, has degree at most~2,
and hence $U(f)-24f-2^{11}f^2$, if non-zero, is a function
$X_0(2)\to\PP^1$ of degree at most~4. Hence one can check that this
function is identically zero by computing the first few terms of its
$q$-expansion and verifying that they are zero. The second identity
follows similarly, or directly from the product formulas:
$$f(\tau) f\Big(\tau + \frac{1}{2}\Big) =
q \prod_{n=1}^{\infty} \frac{1}{(1-q^{2n-1})^{24}}
\times (-q) \prod_{n=1}^{\infty} (1+(-q)^n)^{24}$$
$$=
 - q^2 \prod_{n=1}^{\infty} \frac{(1-q^{2n-1})^{24}
(1 + q^{2n})^{24}}{(1-q^{2n-1})^{24}} = 
- q^2 \prod_{n=1}^{\infty} (1 + q^{2n})^{24} = - f(2 \tau).$$
$\qed$
\end{Proof} 

Using this lemma, we may inductively determine $U(f^k)$ for
all positive $k$. To do this, we observe (by definition) that
$$2 \cdot U(f^k) = f \bigg(\frac{\tau}{2}\bigg)^k
+ f\bigg(\frac{\tau+1}{2}\bigg)^k.$$
Thus, if $X_k := U(f^k)$, then multiplying out, one sees
that the $X_k$ satisfy the
recurrence relation:
$$X_{k+2} - \Bigg(
f\bigg(\frac{\tau}{2}\bigg) +  f\bigg(\frac{\tau+1}{2}\bigg) \Bigg)
 \cdot
X_{k+1} + f\bigg(\frac{\tau}{2}\bigg) f\bigg(\frac{\tau+1}{2} \bigg)
 \cdot X_{k} = 0\ \ (k\geq0).$$
Moreover, From Lemma~\ref{lemma:UF}, we may evaluate
the coefficients of this recurrence to conclude that
$X_0=1$, $X_1=24f+2^{11}f^2$ and, for $k\geq2$,
$$X_k = U(f^k) =  (48f + 2^{12} f^2)  U(f^{k-1}) +
f \cdot U(f^{k-2}).$$
In particular, we note that $U(f^k)$ is a polynomial in $f$ with
integer coefficients and of degree at most~$2k$. These results
are in Emerton's thesis and apparently are originally due to Kolberg.

\begin{df} Define integers $s_{i,j}$, $i,j\in\Z_{\geq0}$, by
$$U(f^j) = \sum_{i=0}^{\infty} s_{i,j} f^i.$$
\end{df}
Note that $s_{i,j}=0$ for $i>2j$. We also note
that $s_{i,j}=0$ for $j>2i$, by comparing the coefficients of $q^i$
in the definition of $s_{i,j}$. We have $s_{0,0}=1$,
$s_{1,1}=24$ and $s_{2,1}=2^{11}$, $s_{1,2}=1$, and $s_{i,j}=0$
in all other cases with $0\leq i\leq 1$ or $0\leq j\leq 1$.
\begin{lemma} \label{lemma:sij} The integers $s_{i,j}$ satisfy the
recurrence relation:
$$s_{i,j} = 48 s_{i-1,j-1} + 2^{12} s_{i-2,j-1} + s_{i-1,j-2}\ \ (i,j\geq2).$$
Moreover, for $i,j\geq1$ and $i\leq2j$, $j\leq2i$ we have an equality:
$$s_{i,j} = \frac{(i+j-1)! 3j \cdot 2^{8i-4j-1}}{(2i-j)! (2j-i)!}.$$
\end{lemma}
\begin{Proof} The recurrence for $s_{i,j}$ follows directly from the
recurrence for $U(f^k)$. The explicit  formula also
satisfies this same recurrence and moreover equals $s_{i,j}$
for $s_{i,1}$, $s_{i,2}$, $s_{1,j}$ and $s_{2,j}$. This suffices
to prove the second equality. $\qed$
\end{Proof}

\medskip

The constant function $1$ is an eigenform for $U$ with
eigenvalue $1$. Thus to determine the spectral theory
of  $U$ on $\Mo$ it suffices to work with the cuspidal
subspace: 
$$\So:=\So(\SL_2(\Z);\Q_2;2^{-1/2})$$
 of $q$-expansions
with zero constant term. In the $p$-adic setting,
\emph{cuspidal}  generalises the notion of vanishing at the
single cusp $\infty$;  thus certain Eisenstein series
(such as the twin form of the usual Eisenstein series
$E_{2k}$)
are considered cuspidal.
Lemma~\ref{lemma:sij} provides us with an explicit
 description of the action of $U$ on $\So$. This allows
us to gain fine control over the spectrum of $U$ on $\So$.

Before we begin the proof of Theorem~\ref{theorem:weight0}, 
we recall some elementary facts about
continuous endomorphisms of Banach spaces. If~$M$ is a Banach
space over $\Q_p$, and $\{e_1,e_2,e_3,\ldots\}$ is a countable subset of~$M$,
then we say that $\{e_1,e_2,e_3,\ldots\}$ is
\emph{an orthonormal Banach basis for $M$} if

\begin{itemize}
\item $|e_i|=1$ for all~$i$,
\item every $m\in M$ can be written uniquely as $m=\sum_{i\geq1}a_ie_i$
for a sequence $a_i\in\Qp$ such that $a_i\to0$ as $i\to\infty$, and
\item if $m,a_i$ are as above, then $|m|=\max_i|a_i|$. 
\end{itemize}

If~$M$ is a Banach space and $\{e_1,\ldots\}$ is an orthonormal
Banach basis for~$M$, and if $\phi:M\to M$ is a continuous $\Q_p$-linear
map, then we define \emph{the matrix of $\phi$} to be the collection
$(c_{i,j})_{i,j\geq1}$ such that $\phi(e_j)=\sum_{j}c_{i,j}e_i$.
The collection $(c_{i,j})$ has the following two properties:

\smallskip

(i) For all $j$, $\lim_{i\to\infty}c_{i,j}=0$

\smallskip

(ii) There exists some $C\in\R$ such that $|c_{i,j}|\leq C$ for all $i,j$.

\smallskip

Conversely, given a collection $(c_{i,j})_{i,j\geq1}$ satisfying
(i) and (ii) above, there is a unique continuous linear map $\phi:M\to M$
with matrix $(c_{i,j})$. Composition of linear maps corresponds to
multiplication of matrices using the usual formula, which one easily
checks to converge because of (i) and (ii) above.

Set $r=1/2$, let~$\Mo$ denote the 2-adic Banach space
of $2^{-r}$-overconvergent $2$-adic modular forms of weight~0 and tame level~1,
equipped with the supremum norm,
and let $\So$ denote its cuspidal subspace.
We prove in the appendix that a Banach basis for~$\Mo$
is $\{1,g,g^2,g^3,\ldots\}$ with $g=2^6f$; we consider $\Mo$ as
being equipped with this basis once and for all. Moreover,
$\So$ has a natural basis given by $\{g,g^2,g^3,\ldots\}$.
We can write the matrix of the operator $U$ on $\So$
as $(u_{i,j})_{i,j\geq1}$, where 
$$u_{i,j} = 2^{6j - 6i} s_{i,j} = 
\frac{(i+j-1)! 3j \cdot 2^{2i+2j-1}}{(2i-j)! (2j-i)!}$$
(and where we interpret this as being zero if $i>2j$ or $j>2i$).
Let $A = (a_{i,j})$, $B= (b_{i,j})$, $D = (d_{i,j})$ ($i,j\geq1$),
be respectively the lower triangular, upper triangular and
diagonal matrices defined as follows:
$$\quad a_{i,j} =
\frac{ 2^{2i-2j}  i!^2 (2 j)!^2 (2j + i - 1)!}
{(2 i)! (i - j)! j! (i + j)! (2j -i)! (3 j-1)!} ,
\ \  2j \ge i \ge j,
\qquad
0 \ \mbox{otherwise.}$$
$$b_{i,j} =
\frac{j}{i} \cdot \frac{2^{2j - 2i}  j!^2 (2 i)!^2 
(2i + j - 1)!}{
(2j)! (j-i)! i! (j+i)! (2i-j)! (3 i-1)!} ,
\ \  2i \ge j \ge i,
\qquad
0 \ \mbox{otherwise.}
$$
and
$$d_{i,i} = \frac{2^{4i + 1} (3i)!^2 i!^2}{3 \cdot (2i)!^4}.$$
Note the symmetry in these formulas. One has
$i \cdot b_{i,j} = j \cdot a_{j,i}$, and $a_{i,i} = b_{i,i} = 1$.

\begin{lemma} The matrices defined by $A$, $B$ and $D$ satisfy (i) above,
and all have coefficients in~$\Z_2$ so also satisfy (ii) above.
Moreover $A \equiv B \equiv Id \mod 2$.
\end{lemma}

\begin{Proof} That $A$, $B$ and $D$ satisfy
(i) is clear because if $i>2j$ then $a_{i,j}=b_{i,j}=d_{i,j}=0$.
Thus it suffices to prove that $A \equiv B \equiv Id \mod 2$ and that
$D$ has coefficients in $\Z_2$.
Recall that $v_2\left(\left(2n\right)!\right)=n+v_2(n!)$ and from this we
see that $v_2(d_{i,i})=1+2v_2\left((3i)!/i!\right)>0$, so $D$ has entries
in $\Z_2$.
Next we deal with $A$. Because $a_{i,i}=1$ for all $i$,
it suffices to prove that the valuation of $a_{i,j}$ is 
positive for $2j\geq i > j\geq1$. We write:
$$a_{i,j} = 6ij \left(\frac{(2j)!}{2^j j!} \right)^2 
\left(\frac{2^i i!}{(2i)!} \right)^2 
\left(\frac{(2i-1)!}{(i+j)!}\right)
\left(\frac{(2j+i-1)!}{(3j)!}\right) \binom{j}{i-j}.$$
Again using that $v_2((2n)!) = n + v_2(n!)$, we see that for $i > j$ the
right hand side is clearly $6$ times a product of terms in
$\Z_2$, and so lies in $2 \Z_2$.

For~$B$ the argument is similar:
since $b_{i,i} = 1$ for all $i$,
it suffices to prove that the valuation of $b_{i,j}$ is
positive for $2i \ge j > i$. We write:
$$b_{i,j} = \frac{j}{i} \cdot a_{j,i} =  6j^2 
\left(\frac{(2i)!}{2^i i!} \right)^2
\left(\frac{2^j j!}{(2j)!} \right)^2
\left(\frac{(2j-1)!}{(i+j)!}\right)
\left(\frac{(2i+j-1)!}{(3i)!}\right) \binom{i}{j-i},$$
and again observe that the right hand side is~6 multiplied
by a product of terms all of which lie in $\Z_2$. $\qed$
\end{Proof}

\medskip
Hence the matrices $A,B,D$ all define continuous endomorphisms
of $\So$, which we also call $A$, $B$ and $D$.

\begin{lemma}We have $ADB = U$. 
\end{lemma}

\begin{Proof} It suffices to show that
$$u_{i,j} = \sum_k a_{i,k} d_{k,k} b_{k,j}.$$
The right hand side of this equation becomes, after
expanding out
and simplifying,
$$ u_{i,j} \cdot
\frac{ 4 \cdot i!^2 j!^2 (2i-j)! (2j-i)!}{(2i)! (2j)! (i+j-1)! }
\sum_k
 \frac{ (2k + i - 1)!
 k (2k + j - 1)!}
{ (i - k)!  (i + k)! (j-k)! (j+k)! (2k-i)! (2k-j)!}.$$
Hence it suffices to prove that
$$\frac{(2i)! (2j)! (i+j-1)!}{ 4 \cdot i!^2 j!^2 (2i-j)! (2j-i)!  } =
\sum_k
 \frac{ (2k + i - 1)!
 k (2k + j - 1)!}
{ (i - k)!  (i + k)! (j-k)! (j+k)! (2k-i)! (2k-j)!}.$$
This identity, however,
follows from classical results; for example, we derive it from
a three term 
specialization of Dougall's ${}_7 F_6$
summation
formula (see \cite{Hardy}, (7.2.3)).
First note that by symmetry, we may assume that $i \le j$.
To force the summation to start at $k = 0$, we let
$k = i - n$. Let $(a)_n = (a)(a+1) \cdots (a+n-1)$.
After repeated application of the formal relations
$$(a)_n:= \frac{(a+n-1)!}{(a-1)!}, \qquad
 (-a)_n = \frac{a! (-1)^n}{(a-n)!}, \qquad
(-a)_n (-a+1/2)_n = \frac{(2a)!}{2^{2n} (2a-2n)!}$$
to transform our sum into hypergeometric form,
the required identity becomes the following:
$${}_7 F_6
\left( \begin{array}{c}
(-i)/2,(1-i)/2,(j-2i)/2,(j-2i+1)/2,1-i,-2i,-i-j \\  
(1-3i)/2,(2-3i)/2,(-2i+1-j)/2,(-2i+2-j)/2,-i,j-i+1 \end{array} ; 1
\right)$$
$$ = \frac{3 \cdot (2i)!^2 (2j)! (i+j-1)! (j-i)! (j+i)!}
{4 \cdot i! j!^2 (2j-i)! (3i)! (2i+j-1)!}.$$
The smallest integer in the
numerator is $i/2$ or $(i-1)/2$, and we consider each
case separately.
Dougall's summation
formula expresses this hypergeometric sum as a
m\'{e}lange
of rising factorials that eventually 
simplify to the required answer. $\qed$
\end{Proof}

\medskip 

An alternative method for proving our identity
would be
via the automated ``creative telescoping''
of Zeilberger.
See for example \cite{Zeil}, where Zeilberger proves Dougall's
summation formula in one (rather long) line; the very short proof
there specialises to a proof
of the identity we require (note however that the statement of the theorem
in~\cite{Zeil} contains a typographical error: the $(-1-a-b-c-d)_n$
term in the denominator should be $(-1-a+b+c+d)_n$ and the proof
should be modified similarly).

\begin{lemma} \label{lemma:NP} The Newton Polygon of $U:=ADB$ is the same as
the Newton Polygon of $D$.
\end{lemma}

\begin{Proof} We use only the fact that $A$ and $B$ are both congruent
to the identity modulo~$2$, and that~$D$ is integral,
diagonal, and compact. Because $B$ is congruent
to the identity mod~2 it has an
inverse. Note  by \S5, Corollaire~2 of~\cite{serre},
$ADB$ has the same Newton Polygon as
$B(ADB) B^{-1} = BAD$, so it suffices to
prove that $CD$ has the same Newton Polygon as $D$, for
any matrix $C$ congruent to the identity modulo $2$.

If $X=(x_{i,j})_{i,j\geq1}$ is
a matrix, and $r_1,r_2,\ldots,r_n$ are distinct positive
integers, then by the $n\times n$ \emph{principal minor} of~$X$ associated
to these integers we mean the determinant of the $n$ by $n$ matrix
formed from the $r_i$th rows and columns of~$X$, $1\leq i\leq n$.
If~$X$ is the matrix associated to a compact morphism,
then the Newton Polygon of $X$ is the lower convex hull of the
points $(n,\Sigma_n)\in\R^2$, where $\Sigma_n$ is the valuation of the sum
of all $n \times n$ principal minors of $M$. 

Firstly, note that if~$d$ is any $n\times n$ principal minor
of the diagonal matrix~$D$, then $(n,v_2(d))$ lies on or above
the Newton Polygon of~$D$.
Secondly, note that for each $(n,\Sigma_n)$ that
lies at a vertex of the Newton Polygon of~$D$,
there is a unique $n\times n$ minor with maximal valuation.
Both of these facts are easily verifiable using the fact that~$D$
is diagonal. 
Next note that if $r_1,r_2,\ldots,r_n$ are distinct positive
integers then the principal minors of $D$ and $CD$ associated
with these integers have the same valuation, because the $n\times n$
minor of~$CD$ associated to these integers is just the product
of the minor associated to~$C$ (which is a unit) and the minor associated
to~$D$. Hence all principal minors of~$D$ and~$CD$ have the same valuation
and now it is easy to check that this forces the Newton polygons
of~$C$ and~$D$ to be the same.

$\qed$
\end{Proof}

Theorem~\ref{theorem:weight0} follows immediately from the above lemma.

\section{Extensions and Generalizations}

\subsection{What is special about the function $f$?}

There are many ways to parameterise a ($p$-adic) disc,
but the choice of~$f$ that we made led
to the simple formulae which enabled us to prove results about slopes.
One possible reason why this~$f$ was a good choice is that the basis 
defined by 
powers of $f$ behaves well with respect to a certain pairing,
which we define below. Recall
our function $g$ defining an isomorphism of~$X_0(2)$ with the
projective line. Let~$w$ denote the Atkin-Lehner involution
on~$X_0(2)$. 

\begin{lemma}\label{g} For $n\in\Z$ we have $\al g^n=g^{-n}.$
\end{lemma}
\begin{Proof}
It suffices to prove the lemma for $n=1$.
Since $\Delta(-1/\tau) = \tau^{12} \Delta(\tau)$, we see that
$$\al g = g(-1/(2 \tau)) =  \frac{2^6  \Delta(-1/\tau)}{\Delta(-1/(2 \tau))} =
\frac{2^6  \tau^{12} \Delta(\tau)}{2^{12} \tau^{12} \Delta(2 \tau)} =
\frac{1}{g}.$$
$\qed$ \end{Proof}

Let~$X$ denote the rigid affinoid annulus
in $X_0(2)$ defined by $X=\{x\in X_0(2):|g(x)|=1\}$.
This is the width zero annulus ``in the middle'' of the supersingular
annulus in $2^{-6}<|g|<2^6$ in $X_0(2)$, and the Atkin-Lehner involution~$w$
induces an involution $X\to X$. If~$\eta$ is a holomorphic one-form
on~$X$ then we can write $\eta=(\sum_na_ng^n)dg$ and we
define $\int_\infty\eta:=a_{-1}.$
Similarly we can write
$\eta=\left(\sum_nb_n(1/g)^n\right)d(1/g)$ and we define $\int_0\eta:=b_{-1}.$
An easy check shows that $\int_0\eta+\int_\infty\eta=0$.

\begin{df} Let $\ip$ denote the following bilinear form on $\Mo$:
$$\<\alpha,\beta\> =  \int_{\infty} \al \alpha \cdot d \beta.$$
\end{df}
\begin{lemma} The bilinear form $\ip$ is symmetric.
\end{lemma}
\begin{Proof}
Since $\al$ swaps the two cusps, we see that
$$\<\alpha,\beta\> =  \int_{0} \alpha \cdot d(\al \beta)$$
Since
$\int_0 + \int_{\infty} = 0$ we see that
$$\<\alpha,\beta\> =  - \int_{\infty} \alpha \cdot d(\al \beta)$$
In particular,
$$\<\alpha,\beta\> - \<\beta,\alpha\> 
= \int_{\infty} \al \alpha \cdot d \beta +   \beta \cdot d(\al \alpha)
= \int_{\infty} d((\al \alpha) \beta) = 0.$$ 
$\qed$ \end{Proof}
\begin{lemma} The basis $\{g^k\}_{k=0}^{\infty}$ is an
orthogonal basis for $\M_0$, with respect to $\ip$.
\end{lemma}
\begin{Proof}
By lemma~\ref{g} we have $\al g^k = g^{-k}$. Hence
$$\<g^m,g^n\> =  \int_{\infty} g^{-m} \cdot
n g^{n-1} dg =
\int_{\infty} n
 g^{n-m} \cdot  \frac{dg}{g} = n \cdot \delta_{m,n}.$$
$\qed$ 
\end{Proof}

Thus
powers of $g$ behave nicely with respect to this
pairing. On the other hand, this pairing behaves nicely with respect
to $U$:

\begin{theorem} $U$ is self-adjoint with respect to $\ip$.
\end{theorem}

\begin{Proof} It is enough to show this for any pair $g^j$, $g^i$.
We see that
$$\<U g^j,g^i\> = \sum_{k=0}^{\infty} u_{k,j} \<g^k,g^i\> =
u_{i,j} \<g^i,g^i\>.
$$
Now from the above calculation and our explicit evaluation of
$u_{i,j}$, this is equal to
$$\frac{2^{2i+2j-1} 3ij (i+j-1)!}{(2i-j)!(2j-i)!} 
= \<g^j,U g^i\>,$$
since the penultimate expression is symmetric.
$\qed$ \end{Proof} 

It is natural to ask whether this theorem is a special case of
a more general phenomenon.

\subsection{Weights other than Zero}

Our results in this section are unfortunately much more
incomplete. We may relate the action of $U$ in weight
$0$ to the action of $U$ in weight $k$ by ``Coleman's
trick'', namely, a judicious application of the
identity
$U(g V(h)) = h U(g)$. In our case we  again obtain explicit
formulae for matrix entries of~$U$ acting on $\Sk_k$. Take
for example the case $k=-12m$ with $m\geq0$ an integer.
Define
$$h_k = \frac{\Delta(2 \tau)^{m}}{\Delta(\tau)^{2m}} =
\Delta(2 \tau)^{-m} f^{2m}$$
By Coleman's trick we now observe that
$$U(h_k f^j) = U(\Delta(2 \tau)^{-m} f^{j+2m}) =
\Delta(\tau)^{-m} U(f^{j+2m}) = h_k f^{-m} U(f^{j+2m})$$
In particular, with respect to the basis
$\{h_kg,h_kg^2,h_kg^3,\ldots\}$, $U$ in weight $-12m$
can be given explicitly by the matrix
$(2^{-6m}u_{i+m,j+2m})_{i,j\geq1}$. Since for $m \in \mathbb N$, $-12m$ is
dense in $4 \Z_2$, if one could prove Conjecture~\ref{conj:nonclassical}
in the case $k=-12m$ then by a continuity argument one could prove
it for all $k$  congruent to~0 mod~4, and then by using
the theta operator one should be able to deal
with the case $k\equiv2$~mod~4 as well.
However, although we have a factorization of the form $U=ADB$, the
matrices $A$ and~$B$ (and $B A$)
 are unfortunately not integral, and new methods
seem to be required.

\subsection{Primes other than $p = 2$}

Our methods rely strongly on the fact that $X_0(p)$ has genus zero.
Presumably our results can be extended to other primes with
this property, with perhaps a corresponding increase in combinatorial
difficulty. On the other hand, our techniques fail, or at least must
be greatly modified, when the genus of $X_0(p)$ is greater than zero,
even though
similar results seem to be true. For example, the following
conjecture appears to be compatible with the conjectures
in~\cite{buzzard:slopeconjectures}:

\begin{conj}
Let $|\lambda_1| \ge |\lambda_2|
 \ge |\lambda_3| \ge \ldots$
be the non-zero eigenvalues \emph{(}with multiplicities\emph{)}
 of $U$ on $\So(\SL_2(\Z);\Q_{11};11^{-r})$ 
for some rational $r$ with $11/12>r>0$.  Then the slope \emph{(}that
is, the $11$-adic valuation\emph{)}
 of $\lambda_n$ is given by the following
formula:
$$v_{11}(\lambda_{n}) =  
v_{11}
 \left(\frac{\left[\frac{6n+1}{5}\right]! \left[\frac{6n+4}{5}\right]!}
{\left[\frac{n}{5}\right]! \left[\frac{n}{5}\right]!}\right)
+ \sum_{k=1}^4 \left[\frac{n+k}{5}\right]
$$
where $[x]$ denotes the largest integer less than or equal to $x$.
\end{conj}

\section{Appendix: Overconvergent forms at small level.}

For this appendix we work with a general level and a general prime~$p$.
For want of a reference, we explain how to extend the theory of
overconvergent $p$-adic modular forms to cases where the level
structure is too coarse for the resulting moduli problem to be rigid.
The motivation is that there is a growing theory of ``explicit''
computations with $p$-adic modular forms, where typically both the
level and the prime are small---for example, this paper. On the other hand,
at several points in theoretical papers on the subject
the hypothesis is made that the level in which one is working is at least~5,
or that $p\geq3$, for convenience (lifting the Hasse invariant,
for example). 

One problem at small level is that there is, as far as we know, no reference
for the theory of ``rigid analytic stacks'', so we proceed by the usual
low-level ad hoc methods. Here is an overview of the idea: take
the level structure one is interested in, and a sufficiently
small auxiliary Galois
level structure with Galois group~$G$ (for
example a full level~$M$ structure for some large prime~$M$). The
product level structure is sufficiently fine to imply the
existence of a compact curve~$X$ representing a moduli problem
on generalised elliptic curves, and~$X$ is equipped with an action of~$G$.
One can form the quotient curve~$X/G$, which is not in general the solution
to a moduli problem, but which is a coarse
moduli space. One cannot always form an appropriate sheaf $\omega$ on these
quotient curves, because $\omega$ does not always descend from~$X$
(problems at elliptic points, for example), but one can still define modular
forms as $G$-invariant sections of tensor powers of~$\omega$ on~$X$.
Next one has to check that the rigid-analytic subspaces that one is
interested in are all $G$-invariant, which comes down to checking that
they have an intrinsic definition that only depends on the underlying
elliptic curve and not the level structure.

We now formulate everything rigorously. Let $U_N$ denote the compact open
subgroup of $\GL_2(\Zhat)$ consisting
of matrices congruent to the identity modulo~$N$. If $U$ is an arbitrary
compact open subgroup of $\GL_2(\Zhat)$ then we define
the \emph{level} of $U$ to be the smallest positive integer $N$ such
that $U_N\subseteq U$.

Let $\Gamma$ be a compact open subgroup of $\GL_2(\Zhat)$
and let $p$ be a prime that does not divide the level of~$\Gamma$.
We shall
give the main definitions in the theory of overconvergent $p$-adic modular
forms for~$\Gamma$. We say that a level structure~$\Gamma$ is
``sufficiently small'' if it satisfies the following
two conditions:

\medskip

(i) The identity
element is the only element of $\Gamma$ which is of finite order and
is conjugate in $\GL_2(\Af)$ to an element of $\GL_2(\Q)$, and

\smallskip

(ii) If $C$ denotes the surjective $\Zhat$-module homomorphisms
$(\Zhat)^2\to\Zhat$, equipped with its natural right action
of $\GL_2(\Zhat)$, then $-1$ operates without fixed points on
$C/\Gamma$.

\medskip

These properties are used in the following way. Property~(i)
implies that there will be no elliptic points in the associated
moduli space, or, more precisely, that the associated
moduli problem on elliptic curves
is rigid and hence representable (\cite{katzmazur}, appendix
to \S4). Property~(ii) implies that there will be no irregular cusps
on the compactification of the representing object (see~\cite{katzmazur},
section 10.13 and Proposition~10.13.4). We remark that
$\Gamma=U_M$ for any $M\geq3$ satisfies properties~(i) and~(ii).
However, perhaps the most common modular curves that one sees are the
curves $X_0(N)$, corresponding to the compact open subgroup $\Gamma$
of matrices in $\GL_2(\Zhat)$ which are upper triangular mod~$N$,
and property (i) always fails for these~$\Gamma$,
since $-Id \in \Gamma$.

Let us assume initially that $\Gamma$ is sufficiently small.
Because $\Gamma$ satisfies (i), the associated moduli problem
on elliptic curves is representable over $\Zp$, by a smooth affine curve
$Y(\Gamma)$, equipped with a universal elliptic curve
$\pi:\E(\Gamma)\to Y(\Gamma)$. Set
$\omega:=\pi_*\Omega^1_{\E(\Gamma)/Y(\Gamma)}$.  Then $\omega$ is
an invertible sheaf on $Y(\Gamma)$. 
Because $\Gamma$ satisfies (ii), the natural compactification $X(\Gamma)$
of $Y(\Gamma)$ is also the solution to a moduli problem (that of
parameterising generalised elliptic curves with level structure---%
see~\cite{deligne-rapoport}, III.6, although we shall not use this
in what follows). Next note that
the sheaf $\omega$ extends in a natural
way to an invertible sheaf on $X(\Gamma)$ (\cite{katzmazur}, 
Proposition~10.13.14). We define
a classical modular form of level $\Gamma$ and weight $k$, defined
over $\Q_p$, to be a global
section of $\omega^{\otimes k}$ on the generic fibre of~$X(\Gamma)$.

If $R$ is a $\Z_p$-algebra then we denote by $X(\Gamma)_R$ the base change
of $X(\Gamma)$ to~$R$. The special fibre $X(\Gamma)_{\Fp}$ of $X(\Gamma)$
is a smooth proper geometrically connected curve,
and has finitely
many supersingular points. Let~$P$ denote a supersingular point
and say $P$ is defined over the finite field~$\F$. Let~$W$ denote
the Witt vectors of~$\F$, and let $K$ denote the field of fractions of $W$.
If $X(\Gamma)^{\an}$ denotes
the rigid space over~$K$ associated to $X(\Gamma)$ then there is a
reduction map from $X(\Gamma)^{\an}$ to $X(\Gamma)_\F(\Fbar)$, and the
pre-image~$U$ of~$P$ is isomorphic to an open disc. The completed local
ring of $X(\Gamma)_W$ at~$P$ is a $W$-algebra
isomorphic non-canonically to a power series ring $W[[t]]$ in one variable;
let us fix one such isomorphism. Then~$t$ can be thought of as giving an
isomorphism from~$U$ to the rigid analytic open unit disc.
Hence if~$r\in\Q_{\geq0}$ then we can
talk about the open subdisc $\{x:|t(x)|<p^{-r}\}$ of~$U$. These subdiscs
in general depend on the fixed isomorphism between the completed local
ring and the power series ring, but if $r<1$ then an easy calculation
shows that they are independent of such choices---the point that we
have chosen to be the centre of~$U$ is a~$K$-point
and the other $K$-points in~$U$ all have distance at most $1/p$ from
our chosen centre, because $K$ is an unramified extension of $\Q_p$.
See~\S3 of~\cite{buzzard:wild} for more details of this construction,
or~\S2 of~\cite{buzzard-taylor}).

The universal formal deformation of the elliptic curve $E_0/\F$
corresponding to~$P$ is an elliptic curve over the completed local
ring of $X(\Gamma)_W$ at~$P$, and hence can be regarded via our
fixed isomorphism as an elliptic curve $E/W[[t]]$. Fix a basis
$\eta$ of $H^0(E,\Omega^1_{E/W[[t]]})$, so
$$H^0(E,\Omega^1_{E/W[[t]]})=W[[t]] \cdot \eta.$$
The Hasse invariant can be thought of as a mod~$p$ section
of the $(p-1)$st tensor power of this module, and hence an element
$A(t)\eta^{\otimes(p-1)}$
of $\F[[t]] \cdot \eta^{\otimes(p-1)}$. It is well-known
that
the Hasse invariant has a simple zero at every supersingular elliptic
curve, which translates into the fact that $A(t)$ is a uniformiser
of $\F[[t]]$. This provides the bridge between our point of view and
Katz'. For example, Katz' analysis of the $p$-divisible group
associated to an elliptic curve and its relation to a lifting
of the Hasse invariant in \S3.7 of~\cite{katz} shows that
for a $\Kbar$-point $u\in U$ corresponding to an elliptic curve $E_u$,
if $|t(u)|>1/p$ then $|t(u)|$ is independent of all choices we have
made, and depends only on the isomorphism class of $E_u$ (and not on the
level structure).

Assume from now on that $0\leq r<1$, and 
that~$\F$ is sufficiently large so that all the
supersingular points in $X(\Gamma)_{\Fp}$ are defined over~$\F$.
Choose parameters $t$ as above for every supersingular point, and
define $X(\Gamma)_{K,\geq p^{-r}}$ to be the rigid space over~$K$ which
is the complement of the open discs $\{x:|t(x)|<p^{-r}\}$ as above,
as $P$ ranges over all supersingular points of $X(\Gamma)_{\F}$.
If $u$ is a $\Kbar$-valued point of $Y(\Gamma)_K$ then let $E_u$
denote the fibre of the universal elliptic curve above $x$. So
$E_x$ is an elliptic curve defined over $K(x)$. Now Katz' arguments
show that if $E_x$ has good supersingular reduction
and if $|t(x)|>p^{-1}$, then $|t(x)|=|t(\sigma x)|$ for any
$\Qp$-automorphism of the field $\overline{K(x)}$. Hence for $r\in\Q$
with $0\leq r<1$ the rigid space $X(\Gamma)_{K,\geq p^{-r}}$ is
the base extension to~$K$ of a rigid subspace $X(\Gamma)_{\geq p^{-r}}$
of $X(\Gamma)$ defined over $\Q_p$. 

Say $r\in\Q$ with $0\leq r<1$.
We define a $p^{-r}$-overconvergent modular
form of level $\Gamma$ and weight~$k$, defined over~$\Qp$,
is a section of $\omega^{\otimes k}$ on $X(\Gamma)_{\geq p^{-r}}$.

Now let $\Gamma$ be an arbitrary (not necessarily sufficiently small)
compact open subgroup of $\GL_2(\Zhat)$,
and let~$p$ be a prime not dividing the level of~$\Gamma$. Choose a
prime $M>2$ dividing neither~$p$ nor the level of~$\Gamma$,
and let $\Gamma'$ denote
$\Gamma\cap U_M$. Then $\Gamma'$ is sufficiently small in the
sense above, and so all of the definitions above apply to $\Gamma'$.
Furthermore, $\Gamma'$ is normal in $\Gamma$; let~$G$ denote
the quotient group. Then~$G$ is finite (in fact $G$ is isomorphic
to $\GL_2(\Z/M\Z)$) and~$G$ acts on $X(\Gamma)$ and $\omega$. Moreover,
because $|t(x)|$ (notation as above) only depends on the elliptic
curve~$E_x$ and not any level structure, $X(\Gamma')_{\geq p^{-r}}$
is $G$-invariant if $0\leq r<1$.

This motivates the following definitions. We define $X(\Gamma)$
to be quotient of $X(\Gamma')$ by the finite group~$G$; note
that $X(\Gamma')$ is a projective curve so taking quotients
is not a problem. In practice,
one can form the quotient in the following manner: $Y(\Gamma')$ is affine,
and is hence of the form $\Spec(R)$, $R$ a ring with an action of~$G$.
One defines $Y(\Gamma)=\Spec(R^G)$ and then compactifies. Note
that the sheaf $\omega$ will probably not in general descend to $X(\Gamma)$.
However, we can still define an $r$-overconvergent modular
form of level $\Gamma$ as being a $G$-invariant element of
$H^0(X(\Gamma')_{\geq p^{-r}},\omega^{\otimes k})$, if $0\leq r<1$.
Sometimes these spaces are zero for trivial reasons---for example
if $\Gamma$ contains $-1$ and if $k$ is odd. If they are not zero,
then they are always infinite-dimensional.

One needs to check that these definitions are independent of the auxiliary
choice of $\Gamma'$. This is not too difficult---if $\Gamma_1'$ and
$\Gamma'_2$ are two choices for $\Gamma'$ with Galois groups $G_1$
and $G_2$ then one sets $\Gamma':=\Gamma'_1\cap\Gamma'_2$ and checks
that both definitions for $X(\Gamma)$
coincide with the quotient of $X(\Gamma')$ by $\Gamma/\Gamma'$, and so on.

One may check without too much difficulty that the
weight~0 $r$-overconvergent forms of level $\Gamma$ are just
the functions on $X(\Gamma)_{\geq p^{-r}}$; this comes from the
fact that one can form quotients of affinoids by finite groups
by looking at invariants, and compatibility of this with the
analytification functor.

One useful result is that if $\Gamma\subseteq\Delta$ both have level prime
to~$p$ then the pre-image of $X(\Delta)_{\geq p^{-r}}$ under the canonical
(forgetful functor) map from $X(\Gamma)$ to $X(\Delta)$ is
$X(\Gamma)_{\geq p^{-r}}$. Only a little harder is the fact that
if $\gamma\Gamma\gamma^{-1}\subseteq\Delta$ and $\Gamma$, $\Delta$
have level prime to~$p$, and $\gamma_p=1$, then the same is true
for the map $X(\Gamma)\to X(\Delta)$ induced by $\gamma$. This is
because $|t(x)|$ depends only on the underlying $p$-divisible
group of the elliptic curve.

The above approach is good for theoretical purposes, but is too abstract
in general to be of much computational use. We now show how to use
these ideas to make the claims of this paper rigorous, thus turning
the argument in this paper from a formal one into a rigorous one.
We set $\Gamma=\GL_2(\Zhat)$ and $p=2$, and write $X_0(1)$ for
$X(\Gamma)$. We now explicitly evaluate $X_0(1)_{\geq p^{-r}}$
for $0\leq r<3/4$. 

It is well-known that the $j$-invariant gives an isomorphism
$X_0(1)\to\PP^1$ defined over $\Z_2$. Let us set $M=3$ (notation
as above), so $\Gamma'=U_3$ and $Y(3):=Y(\Gamma')$
is the modular curve over $\Z_2$ parameterising elliptic curves
equipped with two points of order~3 generating the 3-torsion of
the curve. The generic fibre of this curve is not geometrically
connected but this does not matter. The classical theta
series $\theta:=\sum_{a,b\in\Z}q^{a^2+ab+b^2}=1+6q+6q^3+\ldots$ is a modular
form of
level $\Gamma'$ (in fact it has level $\Gamma_1(3)$) and by the
$q$-expansion principle it is a lift of the mod~2 Hasse invariant.
The special fibre of $X(3):=X(\Gamma')$ has two geometric fibres, both of
genus~0, both defined over $\F:=\F_4$, and both with one $\F$-valued
supersingular point. Let~$P$ denote one of these supersingular points,
and define $W$, $K$, $U$ as above.
Choose an isomorphism~$\iota$ of the complete local ring of $X(3)_W$
at~$P$ with $W[[t]]$,
and let $E/W[[t]]$ denote the universal formal deformation
of the elliptic curve corresponding to the point~$P$. Fix a basis
$\eta$ of $H^0(E,\Omega^1_{E/W[[t]]})$, so
$$H^0(E,\Omega^1_{E/W[[t]]})=W[[t]] \cdot \eta.$$
Now if~$f$ is any modular form over $W(\Fbar)$
of level $\Gamma'$ and any weight~$k$, then $f(E)=h\eta^{\otimes k}$
for some $h\in W[[t]]$. We now think of $\eta$ as being fixed, and
identify a modular form~$f$ with the corresponding function~$h$ as above.
If $u\in U$ then we define $|f(u)|=|h(u)|$
and note that this is independent of the choice of $\eta$.
Moreover, if $|f(u)|>1/2$ then this value is also independent of the
choice of $\iota$. 

Via this fixed isomorphism, $\theta$ can be regarded as an element of
$W[[t]]$ and because $\theta$ lifts the Hasse invariant, we know
that $\theta$
mod~2 in $\F[[t]]$ will be of the form $ut+O(t^2)$ with $u\not=0$.
Now consider the classical level one Eisenstein series
$E_4=1+240(\sum_{n\geq1}\sigma_3(n)q^n)$. Note that the $q$-expansion
of $\theta$ is congruent to~1 mod~2, and hence
$\theta^4$ has $q$-expansion congruent to~1
mod~8. In particular $\theta^4\equiv E_4$ mod~8. This congruence can
be thought be thought of as a congruence of elements of $W[[t]]$.
Now for $u\in U$ with $|E_4(u)|>1/8$ we see that $|\theta(u)|^4=|E_4(u)|$
and hence $|\theta(u)|>2^{-3/4}>2^{-1}$. So
$|t(u)|=|\theta(u)|=|E_4(u)|^{1/4}$.
Conversely, if $|t(u)|>2^{-3/4}$ then $|E_4(u)|=|t(u)|^4>1/8$. We conclude
that if $0\leq r<3/4$ then $X(3)_{\geq 2^{-r}}$ is the subregion
of $X(3)$ where $|E_4|>2^{-4r}$. Moreover, because $j=(E_4)^3/\Delta$,
and $|\Delta(u)|=1$ for all $u\in U$ (as $u$ corresponds to an
elliptic curve with good reduction), we see that for $0\leq r<3/4$,
we have that $X(3)_{\geq 2^{-r}}$ is the region defined by $|j|\geq 2^{-12r}$.
We have proved

\begin{proposition} If $p=2$ and $0\leq r<3/4$ then $X_0(1)_{\geq p^{-r}}$
is the subdisc of the $j$-line defined by $|j|\geq 2^{-12r}$.
\end{proposition}

We now let $g$ denote the modular function $2^6\Delta(2z)/\Delta(z)$.
Recall that this is an isomorphism $X_0(2)_{\Q_2}\to\PP^1$. 
Now one checks that $64/j=g/(4g+1)^3$ and hence the map
$X_0(2)_{\Q_2}\to X_0(1)_{\Q_2}$ induced by the forgetful functor
sends the region $\{|g|\leq 1\}$ to the region $\{|j|\geq|64|\}$.
Moreover, this map preserves $q$-expansions, and
induces an isomorphism between these two regions (one can write
down an inverse, for example, to see this).
We deduce that the disc $\{|g|\leq 1\}$ in $X_0(2)$ is
isomorphic to $X_0(1)_{\geq 2^{-1/2}}$, and that the functions
$\{1,g,g^2,g^3,\ldots\}$ are a Banach basis for $2^{-1/2}$-overconvergent
level~1 weight~0 modular forms. Furthermore, because $g$ vanishes
at infinity, the functions $\{g,g^2,g^3,\ldots\}$ form a basis
for the space of $2^{-1/2}$-overconvergent tame level~1 cusp forms
of weight~0.

Finally, we say a word about other weights $k\equiv0$~mod~12.
Write $k=-12m$ and define $h_k:=\Delta(q^2)^m/\Delta(q)^{2m}$.
Then $h_k$ is a meromorphic section of $\omega^{\otimes k}$
on $X_0(2)$ and a computation of $q$-expansions shows that it
is non-vanishing at infinity. It hence defines a section of
$\omega^{\otimes k}$ on $X(3)_{\geq 2^{-1/2}}$ which is
$G:=\GL_2(\F_3)$-invariant and non-vanishing, and hence a trivialisation
of $\omega^{\otimes k}$ on this rigid space. It is now easy to
show that the $G$-invariant sections of $\omega^{\otimes k}$
are exactly the sections of the form $h_k.s$ for $s$ a function
on $X_0(1)_{\geq 2^{-1/2}}$ and this establishes all the claims
about explicit bases of overconvergent modular forms in this paper.

\end{document}